\documentclass[reqno,11pt]{amsart}
\usepackage{a4wide}
\usepackage{amssymb}
\usepackage{amsmath}
\usepackage[pctex32]{graphics}
\usepackage{graphicx}
\usepackage{bm}

\newtheorem{theorem}{Theorem}
\newtheorem{corollary}[theorem]{Corollary}

\theoremstyle{definition}

\newtheorem{remark}[theorem]{Remark}

\newcommand{\be}{\begin{equation}}
\newcommand{\bd}{\begin{displaymath}}
\newcommand{\ed}{\end{displaymath}}
\newcommand{\ba}{\begin{eqnarray}}
\newcommand{\ea}{\end{eqnarray}}

\def\rr{\mathbb R}

\begin{document}

\title[Nonlinear diffusion: Geodesic Convexity $\Leftrightarrow$ Wasserstein Contraction]{Nonlinear diffusion: Geodesic Convexity is equivalent to Wasserstein Contraction}

\author{Fran\c cois Bolley}
\address{Ceremade, Umr Cnrs 7534, Universit\'e Paris-Dauphine, Place du Mar\'echal de Lattre de Tassigny, 75775 Paris cedex 16, France.} 
\email{bolley@ceremade.dauphine.fr}

\author{Jos\'e A. Carrillo}

\address{Department of Mathematics, Imperial College London, SW7 2AZ - London, United Kingdom.}
\email{carrillo@imperial.ac.uk}


\maketitle


\begin{abstract}
It is well known that nonlinear diffusion equations can be interpreted as a gradient flow in the space of probability measures equipped with the Euclidean Wasserstein distance. Under suitable convexity conditions on the nonlinearity, due to R.~J.~Mc\nobreak Cann~\cite{mccann-advances}, the associated entropy is geodesically convex, which implies a contraction type property between all solutions with respect to this distance. In this note, we give a simple straightforward proof of the equivalence between this contraction type property and this convexity condition, without even resorting to the entropy and the gradient flow structure.

\end{abstract}

\footnote{Date : February 5th, 2014. Mathematics Subject Classification : 28A33, 35B35, 35K10, 37L15, 76S05. Keywords : nonlinear diffusions, optimal transport, entropy, convexity.}

\medskip

We consider the nonlinear diffusion equation
\begin{equation}\label{eq}
\frac{\partial u_t}{\partial t} = \Delta f(u_t), \qquad t >0, \, x \in \rr^d
\end{equation}
where $f(r)$ is an increasing continuous function on $r\in[0, + \infty)$ and $C^2$ smooth for $r>0$ such that $f(0)=0$. The well-posedness theory in $L^1(\rr^d)$ for this equation is a classical matter in the nonlinear parabolic PDEs theory developed in the last 40 years, see \cite{Vaz} and the references therein. Following \cite[Chap. 9]{Vaz}, by a solution we mean a map $u = (u_t)_{t \geq 0} \in C([0,\infty),L^1(\rr^d))$, with $u_t \geq 0$ and mass $\int u_t (x) dx = M$ for all $t$, such that 
\begin{enumerate}
\item[i)] for all $T>0$, the function $\nabla (f \circ u_t) \in L^2((0,T)\times\rr^d)$,
\item[ii)] $u$ weakly satisfies equation \eqref{eq}, i.e., it satisfies the identity
$$
\int_0^\infty\int_{\rr^d} \left\{\nabla (f \circ u_t)\cdot \nabla \chi - u_t \frac{\partial \chi}{\partial t}\right\}\,dx\,dt = 
\int_{\rr^d} u_0(x) \,\chi(0,x)\,dx\,,
$$
for all compactly supported functions $\chi\in C^1([0,\infty)\times\rr^d)$.
\end{enumerate}
For all nonnegative $u_0 \in L^1(\rr^d) \cap L^{\infty}(\rr^d)$ with $\int u_0 (x) dx =M$, there exists a unique such solution to \eqref{eq} with initial datum $u_0$,  see \cite[Sect. 9.8]{Vaz}. Let us normalize the mass $M$ to unity in the rest of the introduction for convenience. 


Equation \eqref{eq} admits the map 
$$
\mathcal U(u) = \int_{\rr^d} U(u(x)) \, dx
$$
as a Liapunov functional. Here the map $U \in C([0, + \infty)) \cap C^3((0, + \infty))$ is defined in a unique way by the relations $f( r) = r U'( r)  - U( r)$  on $(0, + \infty)$ and $U(0)=U'(1)=0$. The map $f$ being increasing is equivalent to $U$ being strictly convex, since $f'( r) = r U''(r )$, and $f$ being positive on $(0, + \infty)$ is equivalent to $\psi : r \mapsto r^d U(r^{-d})$ being decreasing, since $\psi' ( r) = - d r^{d-1} f(r^{-d})$. The map $\mathcal U$ can be extended to the set $\mathcal{P}_{2}(\rr^d)$ of Borel probability measures on $\rr^d$ with finite second moment, by setting $+ \infty$ for non absolutely continuous measures with respect to Lebesgue measure.  The seminal work of R.~J.~McCann \cite{mccann-advances} shows that $\mathcal U$ is (geodesically) displacement convex on the space $(\mathcal{P}_2(\rr^d), W_2)$ if and only if $\psi$ is moreover convex on $(0, + \infty)$, that is, if and only if
\begin{equation}\label{conditionprime}
(d-1) f( r) \leq d \, r \, f'( r), \qquad r>0,
\end{equation}
or equivalently,
\begin{equation}\label{condition}
r \mapsto r^{-1+1/d} \, f( r) \quad \textrm{is \; nondecreasing \; on \;} (0, + \infty)\,.
\end{equation}
We also refer to \cite[Chap. 9]{AGS}, \cite[p. 26]{CJMTU}, \cite[Th. 1.3]{sturm-jmpa} or \cite[Chap. 17]{villani-stflour} for this classical notion. In particular, for the porous medium and fast diffusion equations when $f( r) = r^m$, the condition \eqref{condition} writes $m \geq 1 - 1/d$, which is now classical. Here, $W_2$ is the Wasserstein distance defined for $\mu, \nu$ in $\mathcal{P}_2(\rr^d)$ by
$$
W_2(\mu, \nu)^2 = \inf_{\pi} \iint_{\rr^{2d}} \vert y-x \vert^2 \, d\pi(x,y)
$$
where $\pi$ runs over the set of measures on $\rr^{2d}$ with marginals $\mu$ and $\nu$, see  \cite{AGS,villani-stflour}.

F. Otto \cite{otto-cpde} has interpreted \eqref{eq} as the gradient flow of $\mathcal U$ in $(\mathcal{P}_2(\rr^d), W_2)$ and, for an $f$ satisfying \eqref{condition}, deduced the contraction type property
\begin{equation}\label{contraction}
W_2(u_t, v_t) \leq W_2 (u_0, v_0), \qquad t \geq 0\,,
\end{equation}
for positive smooth  solutions of the Neumann problem for \eqref{eq} on a smooth open bounded set of $\rr^d$, see \cite[Eq.~(133) and Prop.~1]{otto-cpde}. This point of view has been extended by~L.~Ambrosio, N. Gigli and G. Savar\'e \cite{AGS}, through the deep and very general theory of gradient flows of geodesically convex functionals in metric spaces: in particular, for any $u_0$ in $\mathcal{P}_2(\rr^d)$ there exists a unique so-called gradient flow solution $(S_t u_0)_{t \geq 0}$ to \eqref{eq}-\nobreak \eqref{condition}; these solutions satisfy~\eqref{contraction}, a consequence of the displacement convexity of $\mathcal U$ and the key Evolution Variational Inequality. Moreover, $S_t u_0 = u_t$ if $u_0 \in P_2(\rr^d) \cap L^{\infty}(\rr^d)$ due to the uniqueness of the Cauchy problem for weak solutions in \cite[Chap.~9]{Vaz} and gradient flows in \cite[Sect.~6.4]{AS}. 

K.-T. Sturm and M.~von~Renesse \cite{sturm-vonrenesse} have further proved that for the {\it heat equation} on a Riemannian manifold, the geodesic convexity of the Boltzmann entropy is {\it equivalent} to \eqref{contraction}, and to the Ricci curvature being nonnegative. The equivalence of geodesic convexity of the internal energy to curvature-dimension conditions was extended in \cite{sturm-jmpa} to nonlinear diffusions. As shown in \cite{sturm-jmpa}, they imply the contraction property \eqref{contraction} for nonlinear diffusions. Similar ideas were early used for gradient flows in Aleksandrov spaces, see \cite{LL} and the references therein.
In this work, we give the proof of the reverse implication, contraction \eqref{contraction} for nonlinear diffusions implies~\eqref{condition}, for the first time to our knowledge. 

Let us finally remark that the above normalization of the mass is for convenience. Actually, the contraction property \eqref{contraction} holds for all weak solutions with nonnegative initial data in $L^1\cap L^{\infty}(\rr^d)$ with finite second moment and equal initial mass. For this purpose in the rest we extend the $W_2$ distance to nonnegative finite measures with equal mass.

\smallskip

In this note, we give a simple straightforward proof of the contraction property \eqref{contraction} for such $f$, resorting neither to the gradient flow structure, nor even to the functional $\mathcal U$, as in \cite{AGS}  or \cite{otto-cpde}. We first give the {\it equivalence} between condition \eqref{condition} on the nonlinearity and the Wasserstein contraction type property \eqref{contraction} for positive smooth solutions on the ball $B_R = \{ x \in \rr^d, \vert x \vert < R\}:$
 
\begin{theorem}\label{theo}
Let $f$ be an increasing continuous map on $[0, + \infty)$, which is $C^2$ on $(0, + \infty)$ and satisfies $f(0)=0$ and let $R>0$ be fixed. Then McCann's displacement convexity condition \eqref{condition} holds if and only if \eqref{contraction} holds for all solutions $u$ and $v$ to the Neumann~problem
\begin{equation}\label{eqBR}
\frac{\partial u_t}{\partial t} = \Delta f(u_t), \quad t>0, \, x \in B_R, \qquad \frac{\partial u_t}{\partial \eta} = 0, \quad t>0, \, \vert x \vert =R
\end{equation}
with $C^2$ positive initial data $u_0$ and $v_0$ on $\{ \vert x \vert \leq R\}$ respectively, with equal mass.
\end{theorem}

\medskip

Here, $\tfrac{\partial u_t}{\partial \eta}$ denotes the normal derivative of $u_t$ on the sphere $\vert x \vert = R$. We notice that the well-posedness theory for the Neumann problem for nonlinear diffusions is a classical question, see \cite[Chap.~11]{Vaz}. Moreover, the solutions corresponding to positive initial data remain positive and bounded below for all times by a straightforward use of the maximum principle \cite[Th.~11.2]{Vaz}, and thus the solutions we deal with in Theorem \ref{theo} are classical and $C^\infty$ smooth for all times.

\

As a direct consequence of the approximation procedures developed in \cite[Sect.~5.5]{otto-cpde} and \cite[Rem.~18, Prop.~13]{CJMTU}, one can derive the $W_2$ contraction property \eqref{contraction} for solutions of the Cauchy problem:

\

\begin{corollary}\label{coro}
Let $f$ be as in Theorem \ref{theo} and satisfy McCann's displacement convexity condition \eqref{condition}. Then the contraction property \eqref{contraction} holds for all solutions $u$ and $v$ to \eqref{eq} with initial data $u_0$ and $v_0$ in $P_2(\rr^d) \cap L^{\infty}(\rr^d)$ respectively.
\end{corollary}

\

Of course, proving equivalence directly with the contraction property for gradient flows in $\rr^d$ would be more natural. The arguments below formally work on $\rr^d$, however, we lack enough regularity of the solutions to make them rigorous. The rest of this paper is devoted to the proof of the main result in Theorem \ref{theo}.

\

{\bf Proof of Theorem \ref{theo}.}

We first prove the {\it sufficient} part by adapting the strategy of~\cite{bgg-12}. Let $u_0, v_0$ be two $C^2$ positive initial data on $\bar{B}_R = \{ \vert x \vert \leq R\}$ with equal mass. Then so are the solutions $u$ and $v$ of \eqref{eqBR} at all times, by the comparison principle \cite[Th.~11.2]{Vaz}.  For all $t \geq 0$, we let $\xi_t[u] = - \nabla (f \circ u_t) / u_t$ be the velocity field governing the evolution of $u$, written as the continuity equation
$$
\frac{\partial u_t}{\partial t} + \nabla \cdot (u_t \, \xi_t[u]) = 0, \qquad t>0, \, x \in B_R.
$$
We also let $\varphi_t$ be a convex map on $\bar{B}_R$ given by the Brenier Theorem such that $\nabla \varphi_t$ pushes forward $u_t$ onto $v_t$, denoted by $v_t=\nabla\varphi_t \# u_t$, and $\pi_t = (I_{\rr^d}\times \nabla \varphi_t) \# u_t$ is optimal in the definition of $W_2(u_t, v_t)$ (see~\cite{brenier},~\cite[Sect. 6.2.3]{AGS}, or~\cite[Th.~10.28]{villani-stflour} for instance). Here, $I_{\rr^d}$ denotes the identity map. We observe that $u_t=\nabla \varphi_t^*\# v_t$ and $\pi_t=(\nabla \varphi_t^* \times I_{\rr^d}) \# v_t$ for the Legendre transform $\varphi_t^*$ of $\varphi_t$. 

\medskip

Then, for all $t \geq 0$,
\begin{equation}\label{dissipW2}
\frac{1}{2} \frac{d}{dt} W_2^2(u_t, v_t)
=
\iint_{B_R \times B_R} (\xi_t[v](y) -\xi_t[u](x)) \cdot (y-x) \, d\pi_t(x,y).
\end{equation}
This follows from steps 2 and 3 in the proof of \cite[Th.~23.9]{villani-stflour} since for all $t \geq 0$ the solutions $u_t$ and $v_t$ are $C^2$ and positive on the compact set $\bar{B}_R$, so that $\xi_t[u]$ and $\xi_t[v]$ are Lipschitz on $\bar{B}_R$. By the definition of the velocity field $\xi_t$ and the image measure property, this is equal to
\begin{multline}\label{dissipW2.2}
 - \iint_{B_R \times B_R} \Big(\frac{1}{v_t(y)} \nabla (f \circ v_t) (y) -  \frac{1}{u_t(x)} \nabla (f \circ u_t) (x)  \Big)  \cdot (y-x) \, d\pi_t(x,y)
\\
= 
- \int_{B_R} \nabla (f \circ v_t) (y) \cdot (y- \nabla \varphi_t^*(y)) \, dy
+ \int_{B_R} \nabla (f \circ u_t) (x) \cdot (\nabla \varphi_t(x)-x) \, dx.
\end{multline}

Now, the probability densities $u_t$ and $v_t$ are $C^{0, \alpha}$ and bounded from below and above by positive constants on the ball $\bar{B}_R$, so the maps $\varphi_t$ and $\varphi_t^*$ are $C^{2, \alpha}$ on $\bar{B}_R$, see {\bf \cite[Th.~ 5.1]{caff-annals}} and \cite[Th.~12.50]{villani-stflour}. Hence, by integration by parts, the term with $u$ in \eqref{dissipW2.2} is
$$
- \int_{B_R} \! \! f \circ u_t (x)\; (\Delta \varphi_t (x)- d) \, dx + \int_{\vert x \vert = R} \! \! f \circ u (x) \, (\nabla \varphi_t (x) -x ) \cdot \frac{x}{\vert x \vert} \, dx \leq - \int_{B_R} \! \! f \circ u_t \; (\Delta \varphi_t - d).
$$
Here we use the fact that $\nabla \varphi_t (x) \in \bar{B}_R$, so that $\nabla \varphi_t (x) -x$ points inwards $\bar{B}_R$ if $\vert x \vert = R$, and then $(\nabla \varphi_t (x) -x ) \cdot x \leq 0$.
Likewise, the term with $v$ in \eqref{dissipW2.2} is
\begin{eqnarray*}
\leq 
- \int_{B_R} f \circ v_t \, (\Delta \varphi_t^* -d) 
&=&
- \int_{B_R} \! \frac{f \big( v_t (\nabla \varphi_t (x))\big)}{v_t (\nabla \varphi_t(x))} \, (\Delta \varphi^* (\nabla \varphi_t(x))-d) \, u_t (x) \, dx
\\
&=&
- \int_{B_R} \det \nabla^2 \varphi_t(x)\,  f\Big(\frac{u_t(x)}{\det \nabla^2 \varphi_t(x)} \Big)  \, \big(\Delta \varphi_t^* (\nabla \varphi_t(x))-d \big) \, dx
\end{eqnarray*}
by the push forward property and the Monge-Amp\`ere equation that holds in the classical sense, see \cite{caff-annals} and \cite[Ex.~11.2]{villani-stflour} for instance. Hence, we deduce
\begin{multline*}
\frac{1}{2d} \frac{d}{dt} W_2^2(u_t, v_t)
\\
\leq - \int_{B_R} \Big[\det \nabla^2 \varphi_t(x)\,  f\Big(\frac{u_t(x)}{\det \nabla^2 \varphi_t(x)} \Big)  \, \big(\frac{\Delta \varphi_t^* (\nabla \varphi_t(x))}{d}-1 \big) + f(u_t(x)) \, ( \frac{\Delta \varphi_t(x)}{d} - 1) \Big] \, dx.
\end{multline*}

\medskip

Notice that for fixed $x$ and $t$ the bracket in the integral is
$$
p^d \, f(r p^{-d}) (S-1) + f( r) (s-1) 
$$
 where $r = u_t(x)$, $p = (\det \nabla^2 \varphi_t(x))^{1/d}, s = \frac{1}{d} \Delta \varphi_t(x)$ and $S =  \frac{1}{d} \Delta \varphi_t^* (\nabla \varphi_t(x))$. Observe now that  $s \geq p$ by the arithmetic-geometric inequality on the positive eigenvalues of the symmetric matrix $\nabla^2 \varphi_t (x)$. Moreover $\nabla \varphi_t^*(\nabla \varphi_t(x)) = x$ so $\nabla^2 \varphi_t^*(\nabla \varphi_t(x)) \, \nabla^2 \varphi_t(x) = \textrm{Id}$ by differentiation in $x$, that is, $\nabla^2 \varphi_t^*(\nabla \varphi_t(x)) = (\nabla^2 \varphi_t(x))^{-1}$; hence also $S \geq p^{-1}$ for this matrix. Finally $f \geq 0$, so the expression above is
 $$
 \geq p^d f(r p^{-d}) (p^{-1} - 1) + f( r) (p-1) = (p-1) (f( r) - p^{d-1} f(r p^{-d})) \geq 0
 $$
by condition \eqref{condition}. This concludes the proof of the sufficient part.
 
 \
 \
  
 We now turn to the proof of the {\it necessary} part. The strategy is to let $u_0$ and $v_0$ be uniform distributions on concentric balls with similar radii and estimate the dissipation
rate of the $W_2^2(u_0,v_0)$ distance as the difference between the radii tends to zero.
More precisely, let us define $u_0$ to be a suitably normalized characteristic function of a ball with radius $a$, and $v_0$ be a small variation of $u_0$, namely the uniform distribution with the same mass on a centered ball with radius $(1+\delta)a$. Then we will apply the contraction property to these two specific initial data and we will let $\delta\to 0$, recovering the condition on $f$ as the first order term in $\delta$. Similar heuristics were used to show the converse for the characterization of the displacement convexity of the entropy, see \cite[Remark 5.18, Exercise 5.22]{villaniAMS}.

We now define all in detail: let $r>0$ and $0 < a < R$ be fixed, and define $u_0$ by $u_0(x) = r$ on the centered ball of radius $a$ and $u_0(x) = 0$ outside. On the other hand, let $0 < \delta <(R-a)/(2 a)$ be fixed and let $v_0$ be the image measure $\nabla \varphi \# u_0$ of $u_0$ by the map $\nabla \varphi$ where 
$$
 \varphi(x) = \int_0^{\vert x \vert} \psi (z) \, dz, \qquad \vert x \vert \leq R\,,
$$
and $\psi$ is the continuous increasing map defined on $[0,R]$ by
$$
  \psi(z) =
  \left\lbrace
    \begin{array}{ll}
      (1+\delta) z \quad
      &
      \textrm{if} \quad 0 \leq z \leq a
      \\
      (z+ (1+2 \delta)a)/2      \quad
      &
        \textrm{if} \quad a \leq z \leq (1+2\delta)a
         \\
      z      
      &
        \textrm{if} \quad (1+2\delta) a \leq z \leq R.
    \end{array}
  \right.
$$ 
Our aim is to write the contraction type property for these initial data $u_0$ and $v_0$ and $t$ close to $0$, and let $\delta$ goes to $0$. However, since $u_0$ is neither positive nor smooth we need to proceed by approximation. The rest of this proof is devoted to show the technical details to make this density argument possible.

\vskip6pt

{\bf Step 1: Regularizing.-}  To take advantage of the exact expression of the dissipation of the Wasserstein distance as given in \eqref{dissipW2} for smooth positive densities we mollify $u_0$ and $\varphi$ (whence $v_0$) in the following way. 

For given $0<\varepsilon < \varepsilon_0$, where $\varepsilon_0$ will depend only on $a$ and $\delta$, we let $\tilde{u}^{\varepsilon}_0$ be a $C^2$ radial function, nonincreasing on each ray, such that $\tilde{u}^{\varepsilon}_0 (x) = r$ if $\vert x \vert \leq a$ and $= \varepsilon$ if $\vert x \vert \geq a + \varepsilon$; then we let $u^{\varepsilon}_0 = M \tilde{u}^{\varepsilon}_0 / \int \tilde{u}^{\varepsilon}_0$, where $ \int \tilde{u}^{\varepsilon}_0$ is greater than $M$ but tends to $M$ as $ \varepsilon$ goes to $0$. Here $M= c_d a^d r$ is the mass of $u_0$, with $c_d$ being the volume of the unit ball in $\rr^d$.

We also mollify $\psi$ into a $C^3$ increasing map $\psi^{\varepsilon}$ on $[0,R]$ which is equal to $\psi$ on $[0,a]$, $[(1+\varepsilon)a, (1+2\delta - \varepsilon)a]$, and $[(1+2\delta) a, R]$, satisfies
$$
\Vert \psi^{\varepsilon} - \psi \Vert_{L^{\infty}[0,R]} \leq A \varepsilon
$$
for some $A =A(a,\delta)$ and $(\psi^{\varepsilon})'(z) \in [\Lambda^{-1}, \Lambda]$ for some $\Lambda \geq 1+\delta >0$, uniformly in $z \in [0,R]$ and $\varepsilon \in (0, \varepsilon_0)$. Then the map 
$$
 \varphi^{\varepsilon}(x) = \int_0^{\vert x \vert} \psi^{\varepsilon} (z) \, dz, \qquad \vert x \vert \leq R
 $$
is radial, $C^4$ and strictly convex on $\bar{B}_R$ by composition (in particular at $0$ since there it is given by  $\varphi^{\varepsilon}(x) = (1+\delta) \vert x \vert^2/2$). Moreover
$$
\nabla \varphi^{\varepsilon}(x) = \psi^{\varepsilon} (\vert x \vert) \frac{x}{\vert x \vert}
$$
and
$$
\Vert \nabla \varphi^{\varepsilon} - \nabla \varphi \Vert_{L^{\infty}(B_R)} = \Vert \psi^{\varepsilon} - \psi \Vert_{L^{\infty}[0,R]} \leq A \varepsilon.
$$

It is also classical that for all $x \neq 0$ the eigenvalues of the matrix $\nabla^2 \varphi^{\varepsilon}(x)$ are $(\psi^{\varepsilon})'(\vert x \vert)$ (with multiplicity $1$ and eigenvector $x$) and $\psi^{\varepsilon}(\vert x \vert) / \vert x \vert$ (with multiplicity $d-1$ and eigenvectors the vectors orthogonal to $x$), as can be seen by writing the matrix $\nabla^2 \varphi^{\varepsilon}(x)$. But $(\psi^{\varepsilon})' (\vert x \vert) \in [\Lambda^{-1}, \Lambda]$ and $\psi^{\varepsilon}(\vert x \vert) / \vert x \vert \in [1, 1+ \delta]$ by construction, so all eigenvalues of all matrices $\nabla^2 \varphi^{\varepsilon}(x)$ are in $[\Lambda^{-1}, \Lambda]$, where the constant $\Lambda$ is independent of $\varepsilon$.

We finally let $v^{\varepsilon}_0 = \nabla \varphi^{\varepsilon} \# u_0^{\varepsilon}$.

\vskip9pt

{\bf Step 2: Use of the Dissipation of $W_2$.-} The initial datum $u_0^{\varepsilon}$ is $C^2$ and positive on the ball $\bar{B}_R$ and the matrices $\nabla^2 \varphi^{\varepsilon}(x)$ have positive eigenvalues. In particular, the map $\nabla \varphi^{\varepsilon}$ is a $C^1$ diffeomorphism, and $v_0^{\varepsilon}$ is also $C^2$ and positive on $\bar{B}_R$ by the change of variables 
\begin{equation}\label{MA}
v^{\varepsilon}_0 (y) = \Big(\frac{u_0^{\varepsilon}}{\det \nabla^2 \varphi^{\varepsilon}} \Big) \big( ( \nabla \varphi^{\varepsilon})^{-1}(y) \big).
\end{equation}
Hence, \eqref{dissipW2}-\eqref{dissipW2.2} hold, in particular at $t=0$, and the contraction property \eqref{contraction} for the solutions with respective initial data $u^{\varepsilon}_0$ and $v^{\varepsilon}_0$, and $t$ close to $0$ ensures that
\begin{equation}\label{dissipconverse}
0 \geq \int_{B_R} \nabla (f \circ v_0^{\varepsilon}) (y) \cdot (\nabla \varphi^{\varepsilon,*}(y) - y) dy + \int_{B_R} \nabla (f \circ u_0^{\varepsilon}) (x) \cdot (\nabla \varphi^{\varepsilon}(x) - x) dx:=I_1^{\varepsilon}+I_2^{\varepsilon}
\end{equation}
for fixed $a$ and $\delta$, and $\varepsilon >0$ small enough. Here we have let $\varphi^{\varepsilon,*} := (\varphi^{\varepsilon})^*$.

\vskip9pt

{\bf Step 3: Passing to the limit $\varepsilon\to 0$.-} We need to find the limit of $I_1^{\varepsilon}$ and $I_2^{\varepsilon}$ as $\varepsilon\to 0$. We will use weak-strong duality to pass to the limit. Both terms are treated in the same manner, then we present our arguments for $I_1^{\varepsilon}$. We divide it into two substeps:

{\bf Step 3.1.-} First of all, we claim that the first term in the scalar product in $I_1^{\varepsilon}$, i.e. $\nabla (f \circ v_0^{\varepsilon})$, weakly converges to the measure
$$ 
- \displaystyle f (r (1+\delta)^{-d}) \, \frac{y}{\vert y \vert} \, \mu_{ (1+\delta)a } (dy)
$$ 
with $\mu_\alpha$ denoting the uniform measure (with density $1$) on the $d-1$-dimensional sphere of radius $\alpha>0$. The idea here is that if $F^{\varepsilon}$ is a $C^1$ map on $\rr^d$ which is equal to $A$ on the ball $\{\vert x \vert \leq C\}$, equal to $\varepsilon$ on $\{\vert x \vert \geq C(1+\varepsilon)\}$, and say radially decreasing in between, then $\nabla F$ weakly converges to the measure $- A \, \vert y \vert^{-1} \, y \, \mu_{C} (dy)$ as $\varepsilon$ goes to~$0$.

Let indeed $\zeta$ be a vector test function on $\bar{B}_R$. Observing that by construction $v_0^{\varepsilon}$ is constant outside the annulus $(1+\delta) a \leq \vert y \vert \leq (1+2 \delta) a$, we obtain
\begin{multline*}
 \int_{B_R} \nabla (f \circ v_0^{\varepsilon}) \cdot \zeta =  \int_{(1+\delta) a \leq \vert y \vert  \leq (1+2 \delta) a} \nabla (f \circ v_0^{\varepsilon}) \cdot \zeta =  \int_{\vert y \vert = (1+2 \delta) a} f \circ v_0^{\varepsilon} (y) \: \zeta (y) \cdot \frac{y}{\vert y \vert}
\\
 - \int_{\vert y \vert = (1+\delta) a} \! \! f \circ v_0^{\varepsilon} (y) \; \zeta (y) \cdot \frac{y}{\vert y \vert}
\, -  \int_{(1+\delta) a \leq \vert y \vert \leq (1+\delta)a +\Lambda \varepsilon} \! \! f \circ v_0^{\varepsilon} \; \nabla \cdot \zeta
\, -  \int_{(1+\delta)a+ \Lambda \varepsilon \leq \vert y \vert  \leq (1+2 \delta) a} \! \! f \circ v_0^{\varepsilon} \; \nabla \cdot \zeta
\end{multline*}
by integration by parts. 

\smallskip

The first integral is bounded by
$$
f(\varepsilon)  \int_{\vert y \vert = (1+2 \delta) a}  |\zeta (y) |
$$ 
and $f(\varepsilon)$ tends to $0$ as $\varepsilon\to 0$. Hence it tends to $0$. 
The second integral is equal to
$$
f \left(\frac{Mr}{(1+\delta)^d \int \tilde{u}^{\varepsilon}_0} \right) \int_{\vert y \vert = (1+\delta) a}  \zeta (y) \cdot \frac{y}{\vert y \vert} \quad \mbox{ converging to } \quad 
f \left(\frac{r}{(1+\delta)^d} \right) \int_{\vert y \vert = (1+\delta) a}  \zeta (y) \cdot \frac{y}{\vert y \vert} \,,
$$
as $\varepsilon\to 0$. 

To deal with the third integral, recall that $u_0^\varepsilon \leq r$ and that all eigenvalues of all matrices $\nabla^2 \varphi^{\varepsilon}(x)$ are in $[\Lambda^{-1}, \Lambda]$. Then again by \eqref{MA} the regularized target densities verify $v_0^\varepsilon \leq r \Lambda^d$. Therefore, the third integral tends to $0$ since $v^{\varepsilon}_0$ and $\nabla \cdot \zeta$ are bounded uniformly in $\varepsilon$, and the volume of the domain of integration tends to $0$ as $\varepsilon\to 0$.
 
The fourth integral also tends to $0$ since the domain of integration is bounded and the integrand is bounded by a multiple of $f(\varepsilon \Lambda^d)$: Let indeed 
 $$
 \vert y \vert \geq (1+\delta)a+ \Lambda \varepsilon = \vert \nabla \varphi^{\varepsilon} (a) \vert + \Lambda \varepsilon\,.
 $$
Since all eigenvalues of all matrices $\nabla^2 \varphi^{\varepsilon}(x)$ are in $[\Lambda^{-1},  \Lambda]$, then $\nabla \varphi^{\varepsilon}$ is $\Lambda$-Lipschitz. Thus, $\vert y \vert \geq  \vert \nabla \varphi^{\varepsilon} (a+ \varepsilon) \vert$ or equivalently $\vert (\nabla \varphi^{\varepsilon})^{-1}(y) \vert \geq  a+ \varepsilon$ since $\nabla \varphi^{\varepsilon}$ is moreover radial. Therefore, we conclude from \eqref{MA} and the definition of $\psi^\varepsilon$ that
 $$
 v^{\varepsilon}_0 (y) = \frac{u_0^{\varepsilon} \big( (\nabla \varphi^{\varepsilon})^{-1}(y) \big)}{\det \nabla^2 \varphi^{\varepsilon} \, \big( \nabla \varphi^{\varepsilon})^{-1}(y) \big)} \leq \varepsilon \, \Lambda^{d}\,,   
$$
and the claim is proved. 

\

{\bf Step 3.2.-} We will show now that the second term in the scalar product in $I_1^{\varepsilon}$ converges strongly. Using that $\nabla \varphi^{\varepsilon}$ is a diffeomorphism and that $\nabla \varphi^{\varepsilon,*}$ is also $\Lambda$-Lipschitz, we obtain 
\begin{multline*}
\Vert \nabla \varphi^{\varepsilon,*} (y) - \nabla \varphi^* (y) \Vert_{L^{\infty}(B_R)} 
= 
\Vert x - \nabla \varphi^* \big(\nabla \varphi^{\varepsilon}(x) \big) \Vert_{L^{\infty}(B_R)} 
\\
=
\Vert  \nabla \varphi^*\big( \nabla \varphi(x) \big) - \nabla \varphi^* \big(\nabla \varphi^{\varepsilon}(x) \big) \Vert_{L^{\infty}(B_R)} 
\leq 
\Lambda \, \Vert \nabla \varphi - \nabla \varphi^{\varepsilon} \Vert_{L^{\infty}(B_R)}
\leq \Lambda \, A \, \varepsilon.
\end{multline*}
Hence $\nabla \varphi^{\varepsilon,*}$ strongly converges to $\nabla \varphi^*$ in $L^{\infty}(B_R)$.

\vskip9pt

{\bf Step 4: Conclusion.-}  By Steps 3.1 and 3.2 and the weak-strong duality, it follows that $I_1^\epsilon$ in \eqref{dissipconverse} converges to
$$
- f \left(\frac{r}{(1+\delta)^d} \right) \int_{\vert y \vert = (1+\delta) a}  \big( \nabla \varphi^*(y) - y \big) \cdot \frac{y}{\vert y \vert} =   f \left(\frac{r}{(1+\delta)^d} \right) \, \delta \, a^d (1+\delta)^{d-1} c_d
$$
since $\nabla \varphi^*(y) = (1+\delta)^{-1} y$ if $\vert y \vert = (1+\delta) a$.
Analogously, $I_2^{\varepsilon}$ in \eqref{dissipconverse} converges to $-f( r) \delta \, a^d \, c_d$. Therefore, by collecting terms we conclude from passing to the limit in \eqref{dissipconverse}: $0 \geq I_1^\varepsilon + I_2^\varepsilon$, as $\varepsilon\to 0$ that
$$
   0 \geq \delta \, a^d \, c_d \left[ (1+\delta)^{d-1} \, f \left(\frac{r}{(1+\delta)^d} \right) - f( r) \right]
$$
for any fixed $r$ and $\delta>0$. Letting $\delta$ go to $0$ leads to \eqref{conditionprime} at point $r$ by simple Taylor expansion, hence to the equivalent condition \eqref{condition}. This concludes the proof of Theorem \ref{theo}.
\hfill $\Box$

\
 
 \begin{remark}
 R.~J.~McCann's proof of the displacement convexity of the entropy $\mathcal U$ under condition \eqref{condition} in \cite{mccann-advances} is also based on the arithmetic-geometric inequality, though in a much less trivial manner than here.
\end{remark}

\

\begin{remark}
Such contraction properties, or alternatively rates of convergence to equilibrium for the mean field equation
\begin{equation}\label{inter}
\frac{\partial u_t}{\partial t} = \Delta f(u_t) + \nabla \cdot (u_t \, \nabla W * u_t), \qquad t >0, \, x \in \rr^d
\end{equation}
have been derived in \cite{cmcv-03,cmcv-06,BD}. For $W$ strictly convex, and uniformly at infinity only, as in the physically motivated case of $W(z) = \vert z \vert^3$ on $\rr$, only polynomial rates, or exponential but depending on the initial datum, have been obtained by entropy dissipation techniques. Then, bounding from above the (squared) Wasserstein distance between a solution and the steady state by its dissipation along the evolution, universal exponential rates have been derived in \cite{bgg-13} for $f(u)=u$ (linear diffusion). 

The same strategy can be pursued for a nonlinear diffusion. Assume for instance that $d=1$, $f( r)=r^m$ with $1 < m < 2$ and $W$ is a $C^2$ map on $\rr$ for which for all $R>0$ there exists $k( R) >0$ such that $W''(x) \geq k( R)$ for all  $\vert x \vert \geq R$. Then one can prove that for $x_0\in \rr$, equation \eqref{inter} admits a unique steady state $u_{\infty}$ in $\mathcal{P}_2(\rr)$ with center of mass $x_0$, and a constant $c>0$ such that 
$$
W_2(u_t, u_{\infty}) \leq e^{-ct} W_2(u_0, u_{\infty}), \qquad t \geq 0
$$
for all solutions $(u_t)_{t \geq 0}$ to \eqref{inter} with initial datum $u_0$  in $\mathcal{P}_2(\rr)$ with center of mass $x_0$. 

As in \cite{bgg-13}, the proof is based on considering the dissipation of the Wasserstein distance. Two key ingredients consist on one hand in bounding from below the contribution of the (now nonlinear) diffusion term in the dissipation and on the other hand on a good knowledge of the support of the steady state, and of the behaviour of this steady state near the boundary of its support.
\end{remark}

\

\footnotesize \noindent\textit{Acknowledgments.} The authors wish
to thank Marco Di Francesco and Robert J. McCann for discussion on this question, and the referees for relevant remarks and suggestions. This note was partly written while they were visiting ACMAC in Heraklion and finished while the second author was a visiting professor at MSRI in Berkeley. It is a pleasure for them to thank these institutions for their kind hospitality. The first author acknowledges support from the French ANR-12-BS01-0019 STAB project and the second author acknowledges support from projects MTM2011-27739-C04-02 (Spain), 
2009-SGR-345 from AGAUR-Generalitat de Catalunya, the Royal Society through a Wolfson Research Merit Award, and by the Engineering and Physical Sciences Research Council grant number EP/K008404/1.

\

\bibliographystyle{plain}

\end{document}